\documentclass[20pt]{article}
\usepackage{latex8}
\usepackage{epsfig,amssymb}
\input mssymb.tex

\newtheorem{theorem}{Theorem}

\newtheorem{example}{Example}

\newtheorem{definition}{Definition}

\newtheorem{lemma}{Lemma}
\newtheorem{corollary}{Corollary}

\def\R{\Bbb R}
\def\C{\Bbb C}

\def\andd{\hspace{2pt} $\cdot$ \hspace{2pt}}
\def\spp{\vspace{5pt}  \noindent}

\newcommand{\eqref}[1]{{\rm (\ref{#1})}}

\input mssymb.tex

\usepackage{graphicx}
\begin{document}
\title{Nearly Geodesic Riemannian Cubics in $SO(3)$} 


\author{
\bf Lyle Noakes\\
~ \\ 
School of Mathematics and Statistics\\
The University of Western Australia\\
Nedlands, WA 6009\\
Perth, AUSTRALIA\\ 
lyle@maths.uwa.edu.au}
\date{}
\maketitle
\thispagestyle{empty}


\maketitle
\noindent 
{\bf Abstract:}   {\em Riemannian cubics} are curves in a manifold $M$  
that satisfy a variational condition appropriate for interpolation problems. 
When $M$ is the rotation group $SO(3)$, Riemannian cubics are track-summands of {\em Riemannian cubic splines}, used for motion planning of rigid bodies. Partial integrability results are known for Riemannian cubics, and the asymptotics of Riemannian cubics in $SO(3)$ are reasonably well understood. The mathematical properties and medium-term behaviour of Riemannian cubics in $SO(3)$ are known to be be extremely rich, but there are numerical methods for calculating Riemannian cubic splines in practice. What is missing is an understanding of the short-term behaviour of Riemannian cubics, and it is this that is important for applications. The present paper fills this gap by deriving approximations to nearly geodesic Riemannian cubics in terms of elementary functions. The high quality of these approximations depends on mathematical results that are specific to Riemannian cubics.

\spp

\spp
{\bf Mathematics Subject Classification (2000):}\\
 {\em Primary:\/} 70E17 \andd 34E05 \andd 49K99 \andd 70E18 \andd 53A99 \hspace{10pt} {\em Secondary:\/} 70E60 \andd 34H05 \andd 49S05

\clearpage

%
%
\section{Introduction}\label{newintro}
Suppose that a $C^\infty$ curve $x:\R \rightarrow M$ in a Riemannian manifold $M$ is sampled at times $t_0<t_1<\ldots <t_n$, yielding observations $x(t_i)=x_i$ for $0\leq i\leq n$. Then $x$ is uniformly approximated by a track-sum of minimal geodesic arcs joining successive observations, with $O(\delta ^2)$ error, where $\delta $ is the maximum 
distance between $x_{i-1}$ and $x_i$, for $1\leq i\leq n$. So any $C^\infty$ curve in $M$ is a track-sum of curves that are nearly geodesic. Although $x$ is $C^\infty$, the piecewise-geodesic approximation usually fails to be $C^1$ at junctions.  

\spp
A $C^2$ approximation is given by a {\em natural Riemannian cubic spline}, namely a track-sum of {\em Riemannian cubics}, critical for the mean squared norm of covariant acceleration. If sampling is sufficiently frequent we can restrict attention to Riemannian cubics that are {\em nearly geodesic}. 

\spp
In the special case where $M$ is flat, Riemannian cubics are expressed simply in terms of cubic polynomials.
When $M$ is curved, Riemannian cubics are given by a differential equation that is difficult to solve, even when $M$ is the unit $3$-sphere with the standard metric, or the bi-invariant rotation group $SO(3)$. These cases occur in motion planning for rigid bodies, and so numerical methods are needed to find Riemannian cubic splines \cite{barr}. 

\spp
The asymptotics of  Riemannian cubics are studied in \cite{jmp}, \cite{lylenonnull}, \cite{lyleSIAM}, but little is known about the short-term behaviour of cubics, and the short term behaviour is more relevant for motion planning of rigid bodies. The present paper fills this gap, by deriving approximations to nearly geodesic Riemannan cubics in terms of trigonometric functions and polynomials. 
The new approximations are much more informative than Taylor approximation, either in coordinate charts or ambient space, and capture interesting short term geometry of Riemannian cubics that was previously observed in numerical experiments.

\spp
Before describing the methods and layout of the present paper, we review Riemannian cubics in more detail. 
\section{Riemannian cubics}\label{cubicsec}
For $M$ a finite-dimensional Riemannian manifold, consider the functional 
$$J(x)~:=~\int _{t_0}^{t_1}\langle \nabla _tx^{(1)}(t),\nabla _tx^{(1)}(t)\rangle ~dt$$
defined on $C^{\infty}$ curves  $x:[t_0,t_1]\rightarrow M$, where $x$ and its derivative $x^{(1)}$ are prescribed at $t_0,t_1$. Here $\nabla $ denotes the Levi-Civita covariant derivative defined by the Riemannian metric $\langle ~,~\rangle $. 
A {\em Riemannian cubic} is a critical point $x$ of $J$, in the sense that $x$ satisfies the associated $4$th order Euler-Lagrange equation \cite{lylegreg}  
\begin{equation}\label{el} 
\nabla _t^3x^{(1)}+R_{x(t)}(\nabla _tx^{(1)},x^{(1)})x^{(1)}~=~{\bf 0}
\end{equation}
where $R$ is the Riemannian curvature. The existence of Riemannian cubics satisfying the prescribed conditons is proved in \cite{giambo} when $M$ is a complete Riemannian manifold. 
As seen from (\ref{el}), cubically reparameterised geodesics are Riemannian cubics. However most Riemannian cubics do not arise in this way. 
\begin{definition} Given $\epsilon >0$, a $C^\infty$ curve $x:[t_0,t_1]\rightarrow M$ is ($\epsilon $-){\em nearly geodesic} when, for all $t\in [t_0,t_1]$,  
$$\Vert \nabla _tx^{(1)}(t)\Vert ~<~\epsilon \quad \hbox{and}\quad \Vert \nabla _t^2x^{(1)}(t)\Vert ~<~\epsilon .$$ $\square$
\end{definition}
Riemannian cubics need not be nearly geodesic, but the restriction of any $C^\infty$ curve to a sufficiently small subinterval can be reparameterized to a nearly geodesic curve defined over a fixed interval $[t_0,t_1]$. So nearly geodesic curves are informative about the local geometry of arbitrary $C^\infty$ curves, in particular Riemannian cubics.  Nearly geodesic curves arise naturally  in other ways too.  
\begin{example}\label{ex0} Take $M$ to be the matrix group $SO(3)$ of rotations of Euclidean $3$-space $E^3$. Define $x(t)\in SO(3)$ by taking its columns to be the coordinates at time $t$ of an orthonormal frame fixed relative to some rigid body $B$. If the mass distribution of $B$ is spherically symmetric, and if $B$ moves freely, then $x$ is a geodesic in $SO(3)$ with respect to a bi-invariant Riemannian metric. If, however, $B$ is subject to a $C^1$-uniformly small torque $T$ then $x$ is  nearly geodesic. 

\spp
Let $T$ be unknown, and suppose that the configuration of $B$ and its angular velocity are observed at times $t_0,t_1$. Then a minimiser $x$ of $J$ is an interpolant that minimises the mean-squared torque. Since $T$ is $C^1$-uniformly small, the interpolant $x$ is a nearly geodesic Riemannian cubic. $\square$
\end{example}
We refer to \cite{bloch,lyleSIAM,lyletomaszrev,crouch1,crouch4} for further applications of Riemannian cubics.

\spp
When $M$ is Euclidean $m$-space $E^m$, a Riemannian cubic is precisely a polynomial curve of degree $\leq 3$. Nothing like this can be said for non-flat manifolds $M$, even when $M$ is a space of constant nonzero curvature. The situation for {\em elastic curves} is entirely different. 

\begin{example} An {\em elastic curve} in $M$ is a critical point of the restriction of $J$ to the space of {\em constant-speed} curves $x$, 
with $x$ and $x^{(1)}$ still prescribed at $t_0,t_1$. When $M=E^3$ elastic curves are the Euler elastica, whose curvature and torsion are obtained in terms of elliptic sine function, as in Lecture 1 of {\rm \cite{singer}}, whereas Riemannian cubics in $E^3$ are just cubic polynomial curves. 

\spp
On the other hand, when $M$ is the unit sphere $S^3$ in $E^4$, there are quadrature formulae for elastic curves in terms of elliptic functions {\rm \cite{jurd,oscar,elas}} and Lecture 2 of {\rm \cite{singer}}, but quadrature formulae for Riemannian cubics in $S^3$ are known only for a codimension $3$ subclass {\rm \cite{pauley}}. Elastic curves do not resemble Riemannian cubics, except for the very short term. $\square$ 
\end{example}

\spp
Whereas the long term behaviour of Riemannian cubics  is studied in \cite{jmp}, \cite{lyledual}, \cite{lylenonnull}, their short and medium term behaviour is poorly understood. Yet the short and medium term are more significant in applications, such as interpolation and motion planning for rigid bodies. 
\begin{example} Given $t_0<t_1<\ldots <t_n$ and $x_0,x_1,\ldots ,x_n\in M$ define 
$$J(x)~:=~\int _{t_0}^{t_n}\langle \nabla _tx^{(1)}(t),\nabla _tx^{(1)}(t)\rangle ~dt$$
defined on $C^2$ curves $x:[t_0,t_n]\rightarrow M$ satisfying $x(t_i)=x_i$ for $i=0,1,\ldots ,n$. 
A critical point of $J$ is called a {\em natural cubic spline}. Natural cubic splines are characterised as 
$C^2$ track sums of Riemannian cubics on the intervals $[t_{i-1},t_i]$ for $i=1,2,\ldots ,n$, whose covariant acceleration vanishes at $t_0$ and $t_n$. $\square$  
\end{example}
As we shall see in Example \ref{ex1}, Taylor approximations are of limited value in this context. Much more accurate estimates can be made by exploiting  specific properties of Riemannian cubics $x$, especially when $x$ is nearly geodesic. 

\spp
We focus on  bi-invariant $M=SO(3)$ of rotations of $E^3$, and on {\em Lie quadratics} namely the {\em left Lie reductions} $V:[t_0,t_1]\rightarrow E^3$ of a Riemannian cubics $x$. The Lie quadratic $V_\delta$ of a nearly geodesic Riemannian cubic $x_\delta $ is nearly constant. The variational equations of Lie quadratics are used to   find approximations of $V_\delta $. There is a quadrature formulae of \cite{lyledual} for Riemannian cubics in terms of Lie quadratics, but the approximate Lie quadratics cannot be directly substituted for $V_\delta $ into the formula. Nonetheless, taking a little more care, we obtain an approximation  $\hat x$ to $x_\delta $.  

\spp
In effect, the known structural results for Riemannian cubics in $SO(3)$ are exhausted, before  resorting to a Taylor approximation in what is left. Surprisingly, whereas the reconstruction of Riemannian cubics from Lie quadratics requires a quadrature \cite{lyledual}, the first order approximation $\hat x$ for $x_\delta$ given in Theorem \ref{thm4} is algebraic\footnote{The expression is complicated, but it is difficult to see any way around this. } in trigonometric functions and low degree polynomials. 

\spp
The layout of the paper is as follows. 
\begin{itemize}
\item \S \ref{intro} is an introduction to Lie reductions of Riemannian cubics in bi-invariant Lie groups $G$. 
\item \S \ref{approxsec} studies the variational equation of a Lie quadratic, giving examples where the variational equation can  be solved exactly. Example \ref{rconstex} concerns the simplest case of variations through nearly constant Lie quadratics. Derivatives $V^{(i,0)}$ to order $0\leq i\leq n$ of a variation to order $n$ with respect to the variation parameter, give rise to an order $n$ approximation $\hat V_{n}$ of a nearly constant Lie quadratic $V_\delta $.  
\item In \S \ref{nearconstv},~ ${\cal G}$ is taken to be the Lie algebra $so(3)\cong E^3$ of the rotation group $SO(3)$. Theorems \ref{thm2}, \ref{thm3} give formulae for the first and second order approximations $\hat V_{1}$ and $\hat V_{2}$ to the Lie quadratic $V_\delta$. In applications such as Example \ref{ex0}, $V_\delta$ gives the angular momentum relative to the body $B$. 

\spp
Already $\hat V_{1}$ captures significant medium term behaviour of $V_\delta $, as illustrated in Figures \ref{fig:approx0}, \ref{fig:approx1} of Example \ref{ex1}. Taylor approximations to $V_\delta $ perform badly (Figure \ref{fig:approx0}), while $\hat V_{1}$ and $\hat V_{2}$ are nearly indistinguishable from $V_\delta$ in the medium term. For the longer term, $\hat V_{2}$ significantly improves on $\hat V_{1}$ (Figure \ref{fig:approx1}). 
\item In \S \ref{approxsecv} the approximations $\hat V_n$ for $n=1$ and $n=2$ are 
used in combination with Theorem 5 of \cite{lyledual} to derive a first order approximation $\hat x$ to the Riemannian cubic $x_\delta $, namely $x_\delta ^{(j)}(t)=\hat x^{(j)}(t)+O(\delta ^{2})$ for any non-negative integer $j$. This is our main result, stated as Theorem \ref{thm4}, and implemented in Example \ref{ex2}.\end{itemize}
\section{Riemannian cubics in Bi-Invariant Lie groups: Lie quadratics}\label{intro} Now we take $M$ to be a path-connected finite-dimensional Lie group $G$, with bi-invariant pseudo-Riemannian metric. The restriction of the metric to the Lie algebra ${\cal G}$ is an {\em ${\rm ad}$-invariant} semi-definite inner 
product, namely ${\rm ad}(u):{\cal G}\rightarrow {\cal G}$ is skew-adjoint for all $u\in {\cal G}$. Conversely an ${\rm ad}$-invariant semi-definite inner product\footnote{These do not exist for some Lie groups. However for $G$ semisimple we may use the Killing form.} on ${\cal G}$ extends by left multiplication to a bi-invariant 
pseudo-Riemannian metric on $G$. Define the {\em left Lie reduction} $V$ of $x:[t_0,t_1]\rightarrow G$ by 
$$\displaystyle{V(t):=dL(x(t)^{-1})_{x(t)}(x^{(1)}(t)}$$ 
where $L(g)$ denotes left multiplication by $g\in G$. 
\begin{theorem}[\cite{lylegreg,jmp}]\label{lieeqthm} $x$ is a Riemannian cubic 
in $G$ if and only if, for all $t\in [t_0,t_1]$,  and some $C\in {\cal G}$, we have  
\begin{equation}\label{lieeq}
V^{(2)}(t)~=~[V^{(1)}(t),V(t)]+C.
\end{equation}
$\square$
\end{theorem}
Equation (\ref{lieeq}) is  second order, whereas the Euler-Lagrange equation 
(\ref{el}) has order $4$.  Equivalently we may write
\begin{equation}\label{lieeqp}
V^{(3)}(t)~=~[V^{(2)}(t),V(t)].
\end{equation}
A curve $V:[t_0,t_1]\rightarrow {\cal G}$ satisfying (\ref{lieeq}) for all $t$ is said to be a {\em Lie quadratic}. The Lie quadratic and Riemannian cubic are said to be {\em null} when $C={\bf 0}$, and {\em non-null} otherwise. Null Lie quadratics in $E^3$ appear in applications in fluid dynamics \cite{vega}. 

%
\spp
If $x:[t_0,t_1]\rightarrow G$ is a Riemannian cubic then so are $x^{-1}$ and $gx$ 
where $g\in G$. Geodesics in bi-invariant Lie groups are precisely the Riemannian cubics with constant Lie quadratics.
\begin{example}\label{abelianex} An affine line in ${\cal G}$ is a Lie quadratic, and is null 
when its image contains ${\bf 0}$. 
A Lie algebra ${\cal G}$ is abelian precisely when all its affine lines are null Lie quadratics. Then all nonconstant null Lie quadratics are affine lines, and a Lie quadratic  
is precisely a polynomial curve in ${\cal G}$ of degree at most $2$. $\square$
\end{example}
From Theorem \ref{lieeqthm} a short calculation shows 
\begin{corollary} Let $y,z:[t_0,t_1]\rightarrow G$ be restrictions of $1$-parameter subgroups of $G$, with infinitesimal generators $A,B\in {\cal G}$. Then the pointwise product $t\mapsto y(t)z(t)$ is a Riemannian cubic 
if and only if $[[[B,A],A],B]={\bf 0}$. $\square$
\end{corollary}

\begin{example} Let $G$ be bi-invariant $SO(3)$, with $A,B\in so(3)$ linearly independent. Then the pointwise product $t\mapsto {\rm exp}(tA){\rm exp}(tB)$ is a Riemannian cubic if and only if  
$\langle A,B\rangle =0$. $\square$
\end{example}
%
%
%
Null Lie quadratics in $E^3$ are studied in {\rm \cite{jmp,lyleSIAM}}. 
%
%
The present paper  
focuses on the {\em non-null} Lie quadratics and Riemannian cubics, which are generic, and whose behaviour can be complicated, as illustrated in \cite{lylenonnull}. 
%
%
%
\section{Variational Equations For Lie Quadratics}\label{approxsec}
For $t_0<t_1$ and any $C^{\infty}$ map $V:[-1 ,1]\times [t_0 ,t_1] \rightarrow {\cal G}$, define $V_h:[t_0,t_1]\rightarrow {\cal G}$ by $V_h(t)=V(h,t)$, where $h\in [-1,1]$. Denote the $i$-fold derivative of $V$ with respect to $h$ of the $j$-fold derivative with respect to $t$, by $V^{(i,j)}$. 

\spp
Suppose that, for each $h\in [-1,1]$,  $V_h$ is a Lie quadratic. 
Differentiating ~$\displaystyle{V^{(0,3)}=[V^{(0,2)},V^{(0,0)}]}$~ 
$n\geq 1$ times with respect to $h$, we obtain the third order linear ODE for $Y_n(t):=V^{(n,0)}(h,t)$ 
\begin{equation}\label{linearforvar}Y_n^{(3)}+{\rm ad}(V^{(0,0)})Y_n^{(2)}-{\rm ad}(V^{(0,2)})Y_n~=~\sum_{i=1}^{n-1}\left( \begin{array}{c}n\\
i\end{array}\right) [V^{(n-i,2)},V^{(i,0)}]\end{equation}
where the coefficients and right hand side are in terms of  derivatives with respect to $t$ of $Y_i$ with $i=0,1,2,\ldots ,n-1$.  

\spp
Let $V_0$ be the Lie quadratic of a nontrivial cubically reparameterised geodesic, namely $V_0=q_0(t)D$ for ${\bf 0}\not= D\in {\cal G}$ and $q_0:\R \rightarrow \R$ is a quadratic polynomial. For $n=1$, the variational equation (\ref{linearforvar}) evaluated at $h=0$ becomes 
\begin{equation}\label{weq1}Y_1^{(3)}(t)+{\rm ad}(D)(q_0(t)Y_1^{(2)}(t)-q_0^{(2)}Y_1(t))~=~{\bf 0}.\end{equation}
Let ${\rm ad}(D):{\cal G}\rightarrow {\cal G}$ be diagonalizable\footnote{In a semisimple Lie algebra \cite{vara}, the so-called {\em semisimple} elements $D$ for which this holds comprise an open dense subset of ${\cal G}$.} over $\C$, namely ${\cal G}={\cal K}+{\cal E}$ where ${\cal K}$ is the kernel of ${\rm ad}(D)$, and ${\cal E}$ has a basis of eigenvectors $E_i$ of ${\rm ad}(D)$. 
For $K\in {\cal K}$, (\ref{weq1}) gives 
$$\langle Y_1^{(3)}(t),K\rangle =-\langle {\rm ad}(D)(q_0(t)Y_1^{(2)}(t)-q_0^{(2)}Y_1(t)),K\rangle =
\langle q_0(t)Y_1^{(2)}(t)-q_0^{(2)}Y_1(t), {\rm ad}(D)K\rangle =0.$$
So $\langle Y_1(t),K\rangle =q_K(t)$ where 
$q_K:[t_0,t_1]\rightarrow \R$ is another quadratic polynomial. 
Similarly, $\langle Y_1(t),E_i\rangle$ satisfies 
\begin{equation}\label{eqy}y^{(3)}(t) ~=~r(t)y^{(2)}(t)-r^{(2)}y(t)~\Longrightarrow ~ y^{(2)}(t)~=~r(t)y^{(1)}(t)-r^{(1)}(t)y(t)+c_0\end{equation}
where $r(t):=\lambda _iq_0(t)$, $\lambda _i\in \C$ is the eigenvalue of of $E_i$, and $c_0$ is constant. 
The linear ODE (\ref{eqy}) is solvable by quadratures in terms of solutions of the associated homogeneous ODE 
\begin{equation}\label{eqyh}
y^{(2)}(t)~=~r(t)y^{(1)}(t)-r^{(1)}(t)y(t).
\end{equation}
\begin{example}\label{rconstex} Let $r$ be a constant $0\not= k\in \C$. Then (\ref{eqyh}) reads   
$y^{(2)}(t) =ky^{(1)}(t)~\Longrightarrow ~y(t)=c_2+c_1e^{kt}$
where $c_1,c_2$ are constant. So the general solution of (\ref{eqy}) is 
$$y(t)~=~-\frac{c_0}{k}t+c_2+c_1e^{kt}.$$
$\square$
\end{example}
\begin{example}\label{rlinex} Let $r$ be nonconstant and linear. After a time-shift, write 
$r(t)=kt$ where $0\not= k\in \C$ is constant. Then (\ref{eqyh}) becomes 
$y^{(2)}(t) =kty^{(1)}(t)-ky(t)$ and, 
substituting $y(t)=tz(t)$,  
$$2z^{(1)}+tz^{(2)}~=~kt^2tz^{(1)}~\Longrightarrow ~z^{(1)}(t)~=~c_2t^{-2}e^{kt^2/2}$$
from which $z(t)$ is found by quadrature. $\square$
\end{example}

\begin{example}\label{rquadex} Let $r=kt^2$ where $k$ is constant. Equation (\ref{eqyh}) becomes 
$y^{(2)}(t)=kt^2y^{(1)}(t)-2kty(t)$. Substituting $y(t)=w(s)$ with $s=kt^3/3$, gives a form of 
Kummer's equation 
$$sw''~=~(s-\frac{2}{3})w'+\frac{2}{3}w$$
where $w'$ and $w''$ are derivatives with respect to $s$. So $w$ is given in terms of confluent hypergeometric functions. $\square$  
\end{example}
\spp
Suppose now that $V_0$ is constant, and let $h\mapsto V(h,0)$ be an affine function of uniform 
norm $\leq 1$. From equation (\ref{eqy}), for any non-negative integers $n$ and $j$, for some some $c_{n,j}>0$, and for all $(h,t)\in [-1,1]\times [t_0,t_1]$, 
$$\Vert V^{(n,j)}(h,t)\Vert ~<~c_{n,j}.$$
For $\delta \in (0,1)$ and $n\geq 1$, the {\em $n$th order approximate quadratic} ~$\hat V_n:[t_0,t_1]\rightarrow {\cal G}$ ~of~ $V:[-1,1]\times [t_0,t_1]\rightarrow {\cal G}$~  is defined by 
$$\displaystyle{\hat V_{n}(t)~:=~V_0(t)+\sum_{i=1}^n\frac{\delta ^iV^{(i,0)}(0,t )}{i!}}.$$ 
By Taylor's Theorem, 
$\displaystyle{\hat V_n^{(j)}(t)=V_\delta ^{(j)}(t)+O(\delta ^{n+1})}$
where the asymptotic constants depend on $V_0$, $n$ and $j$. Once
${\cal G}$ is specified the analysis of Example \ref{rconstex} can be taken further. 
\section{Approximating Nearly Constant Lie Quadratics in $so(3)$}\label{nearconstv}
Considering Euclidean $3$-space $E^3$ as a Lie algebra, with respect to the cross-product as Lie-bracket, the Euclidean inner product $\langle ~,~\rangle $ is {\rm ad}-invariant. The Lie isomorphism ${\rm ad}:E^3\rightarrow so(3)$, given by ${\rm ad }(v)(w):=v\times w$, identifies $E^3$ with $so(3)$. In particular, an {\rm ad}-invariant inner product $\langle ~,~\rangle $ on $so(3)$ is defined by requiring ${\rm ad}$ to be an isometry. 
Take $G=SO(3)$ with the corresponding bi-invariant Riemannian metric. 

\spp
We seek $n$th order approximate quadratics $\hat V_{n}:[t_0,t_1]\rightarrow so(3)$ of variations $V$, where $V_0$ is a nonzero constant $D\in so(3)$ and $h\mapsto V(h,0)$ is affine.   
Then ${\rm ad}(D)$ is diagonalizable, the linear span ${\cal H}$ of $D$ is a Cartan subalgebra, and ${\rm ad}(D): so(3)\otimes _{\R} \C\rightarrow so(3)\otimes _{\R} \C$ has eigenvalues $\pm d{\bf i}$, with unit eigenvectors $E_1,E_2=\bar E_1\in so(3)\otimes _{\R }\C$ orthogonal to $D$. Define unit vectors   
$$F_0~:=D/d,\quad F_1~:=~(E_1+E_2)/\sqrt{2},~ \quad 
F_2~:=~-(E_1-E_2)/(\sqrt{2}{\bf i})$$
with $[F_0,F_1]=F_2$ and $[F_2,F_0]=F_1$. Then $F_0,F_1,F_2$ corresponds under ${\rm ad}$ to a positively oriented orthonormal basis of $E^3$. 

\spp
For $n\geq 1$ define $f_n:[t_0,t_1]\rightarrow \R$ and $v_n:[t_0,t_1]\rightarrow F_0^\perp$ by 
$\displaystyle{f_n(t)F_0+v_n(t):=Y_n(t)}$.~ Taking 
$k=\pm d{\bf i}$ in Example \ref{rconstex}, $f_1$ is a quadratic polynomial $q$, and $v_1=A+\tilde e(B)$ where $A:[t_0,t_1]\rightarrow F_0^\perp \subset so(3)$ 
is affine, $B\in so(3)$, and $\tilde e(t)$ is the Lie endomorphism ${\rm exp}(-d(t-t_0){\rm ad}(F_0))$ of $so(3)$. 
We have proved
\begin{theorem}\label{thm2} The first order approximate quadratic of $V$ has the form 
$$\hat V_{1}(t)~=~D+\delta (q(t)F_0+A_0+(t-t_0)A_1+\tilde e(t)(B))$$
where $q(t)=c_0+c_1(t-t_0)+c_2(t-t_0)^2$, $c_0,c_1,c_2,d\in \R$ and $A_0,A_1,B\in F_0^\perp$. Then, for $j\geq 0$ and all $t\in [t_0,t_1]$,    
$$V_\delta ^{(j)} (t)~=~\hat V_1^{(j)}(t) +O(\delta ^2)$$
where the asymptotic constants depend only on $d,t_0,t_1$ and $j$. $\square$
\end{theorem}

\spp
For $n\geq 2$, integration of equation (\ref{linearforvar}) gives 
\begin{eqnarray}
\label{forfn}f_n^{(2)}&=&\sum_{i=1}^{n-1}\left( \begin{array}{c}n\\ i\end{array}\right)I([v_{n-i}^{(2)},v_i])\\
\label{forvn}v_n^{(2)}&=&\sum _{i=1}^{n-1}\left( \begin{array}{c}n\\ i\end{array}\right)\tilde i\tilde eI(f_{n-i}^{(2)}\tilde e^{-1}v_i-f_i\tilde e^{-1}v_{n-i}^{(2)}) 
\end{eqnarray}
where ~$\tilde i:={\rm ad}(F_0)$~ and, for a continuous curve $g$ of linear endomorpisms of $so(3)$,  $I(g)(t):=\int _{t_0}^tg(s)~ds$. \\
Then $\tilde i$ commutes with $\tilde e$, and $\tilde i^2=-{\bf 1}$ where ${\bf 1}$ is the identity on $so(3)$. 
If $g$ is $C^1$, then 
\begin{equation}\label{ibp}
I(\tilde eg)~=~\frac{\tilde i}{d}(\tilde eg-g(t_0)-I(\tilde eg^{(1)})),\quad \hbox{in particular}\quad 
I(\tilde e)~=~\frac{\tilde i}{d}(\tilde e -{\bf 1}).
\end{equation}
Taking $n=2$ in equation (\ref{forfn}), and writing $A(t)=A_0+(t-t_0)A_1$ with $A_0,A_1\in so(3)$, repeated use of (\ref{ibp}) gives 
\begin{eqnarray*}
f_2&=&-2\langle [A_0,L_0(B)]+[A_1,L_1(B)],F_0\rangle \\
v_2&=&2q^{(2)}(M_0(A_0)+M_1(A_1)-M_B(B))+2d^2\tilde i\circ I^2(I(q)\tilde e)(B)
\end{eqnarray*}
where  $L_0(t),L_1(t),M_0(t),M_1(t),M_B(t)$ are the endomorphisms of $so(3)$ given by 
\begin{eqnarray}
L_0&:=& ~\frac{1}{d}(-u{\bf 1} +(u^2/2-1)\tilde i +\tilde i\circ \tilde e)\\
L_1&:=&~ \frac{1}{d^2}((u^2/2-3){\bf 1}+ 2u\tilde i+ 3\tilde e +u\tilde i\circ \tilde e)\\
M_0&:=&~\frac{1}{d^3}((u^2/2-1){\bf 1}+ u\tilde i + \tilde e) \\
M_1&:=&~\frac{1}{d^4}((u^3/6-u){\bf 1} + (u^2/2-1)\tilde i+\tilde i\circ \tilde e)\\
M_B&:=&~\frac{1}{d^3}(2(\tilde e-{\bf 1})+u\tilde i\circ (\tilde e+{\bf 1})) 
\end{eqnarray}
with $\tilde e$ evaluated at $t$, and $u:=d(t-t_0)$. This proves
\begin{theorem}\label{thm3} The second order approximate quadratic of $V$ is  
$$\hat V_{2}(t)~=~D+\delta (q(t)F_0+A_0+(t-t_0)A_1+\tilde e(t)(B))+\frac{\delta ^2}{2}(f_2(t)F_0+v_2(t))$$
where $q,A_0,A_1,B$ are as before. For $j\geq 0$ and all $t\in [t_0,t_1]$,  
$$V_\delta ^{(j)} (t)~=~\hat V_2^{(j)}(t) +O(\delta ^3)$$
where the asymptotic constants depend only on $d,t_0,t_1$ and $j$. $\square$
\end{theorem}
Notice that $q,A,B$ contain $3+4+2$ scalar variables, sufficient for  initial conditions of the $3$rd order ODE (\ref{lieeqp}) in $so(3)$. 
%
%
\begin{figure}[p] 
   \centering
   \includegraphics[width=4in]{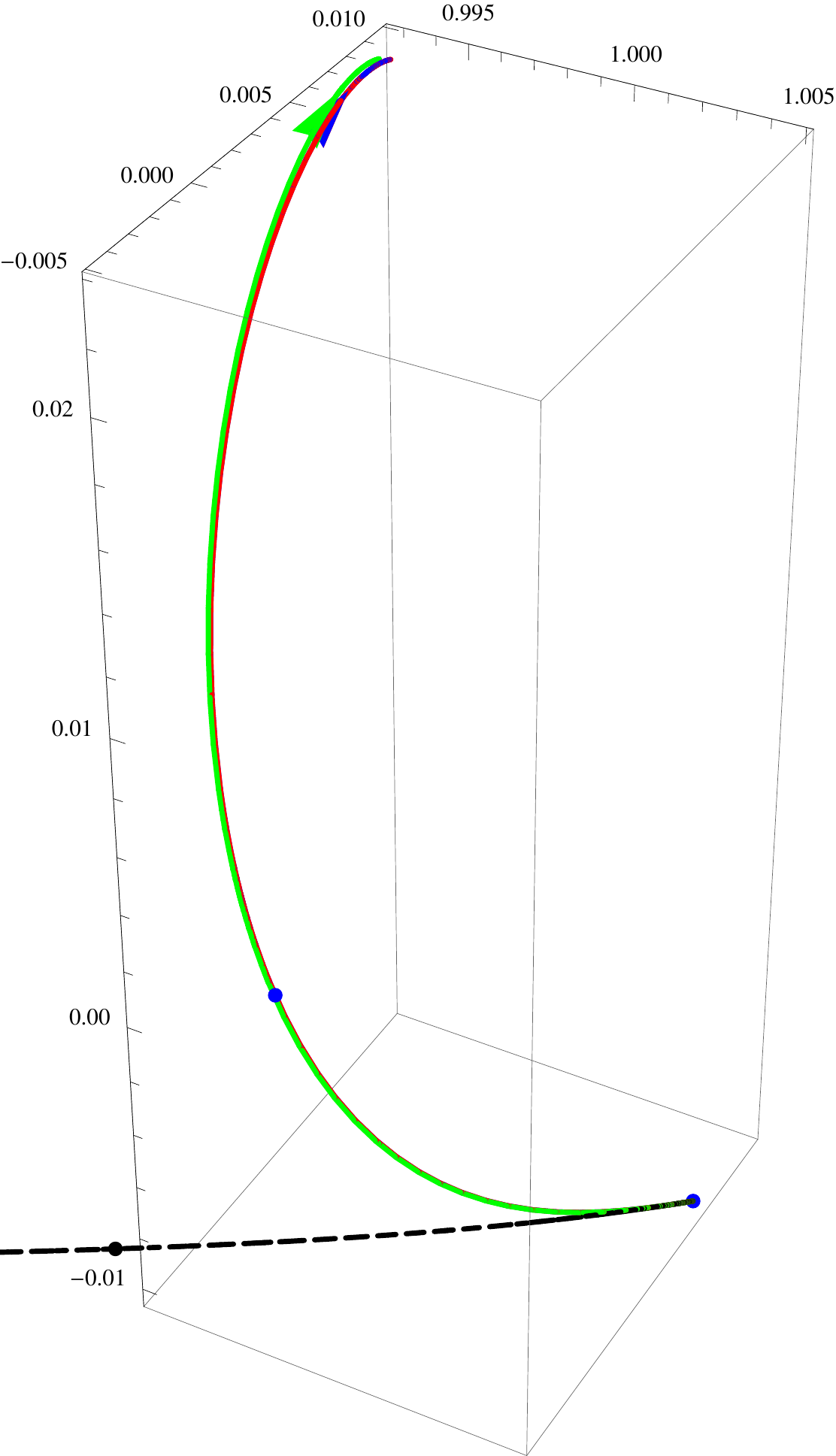} 
   \caption{\rm $V_\delta (t)$ (blue), $\hat V_{1}(t)$ (green), $\hat V_{2}(t)$ (red) and degree $2$ Taylor (dashed)  for $t\in [0,5]$}
   \label{fig:approx0}
\end{figure}

\begin{figure}[htbp] 
   \centering
   \includegraphics[width=7in]{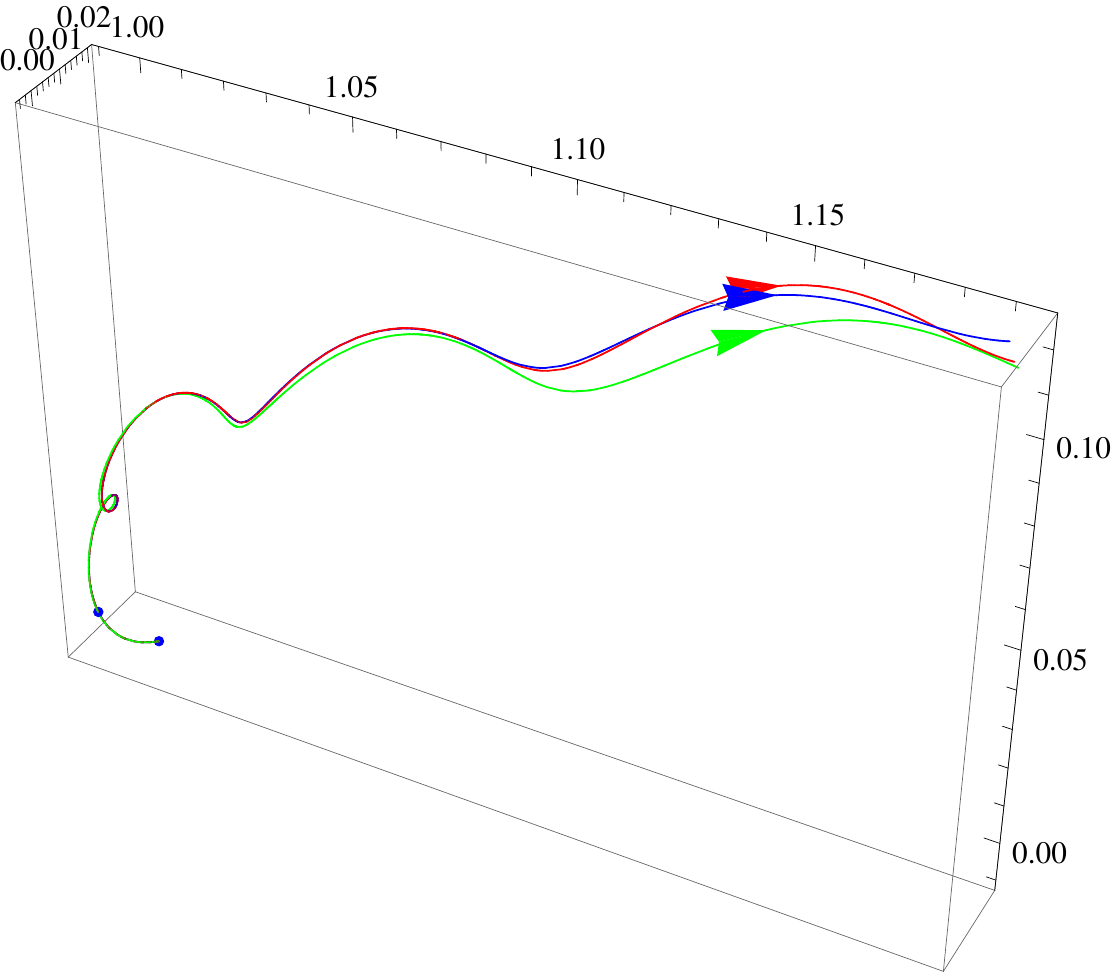} 
   \caption{\rm $V_\delta (t)$ (blue), $\hat V_{1}(t)$ (green), $\hat V_{2}(t)$ (red)  for $t\in [0,25]$}
   \label{fig:approx1}
\end{figure}
\begin{example}\label{ex1} Using ${\rm ad}$ to identify $so(3)$ with $E^3$, Figure \ref{fig:approx0} plots (blue) a numerical solution of (\ref{lieeq}) for a Lie quadratic $V_\delta :[0,5]\rightarrow so(3)\cong E^3$ near $V_0=(1,0,0)$. The numerical solution is obtained using Mathematica's NDSolve, and $V_\delta$ is the non-null Lie quadratic specified by the initial conditions  
$$V_\delta (0)~=~(1.005, 0.006, -0.01),\quad V_\delta ^{(1)}(0)~=~(-0.005, -0.00449, 0),$$
$$V_\delta ^{(2)}(0)~=~(0.001, -0.005, 0.005),\quad V_\delta ^{(3)}(0)~=~(0.00002, 0.005035, 0.005031)$$
with $C\approx (0.0009551, -0.00495, 0.00051755)$. 

\spp
The approximations (green) $\hat V_{1}$ and (red) $\hat V_{2}$ are nearly indistinguishable from $V_\delta$.  The second degree\footnote{Higher order Taylor polynomials are almost as  uncompetitive. } Taylor polynomial, also shown (dashed) in Figure \ref{fig:approx0}, is a poor approximation to $V_\delta $. The initial dot labels $V_\delta (0)$, and the second dots correspond to $t=2$.
%

%
\spp
Figure \ref{fig:approx1} plots the same $V_\delta $ (blue), $\hat V_{1}$ (green) and $\hat V_{2}$ (red) for $t\in [0,25]$.  At first, both approximations successfully follow the contortions of $V_\delta $, but as $t$ increases $\hat V_{1}$ loses accuracy.  The initial dot labels $V_\delta (0)$, the second dot corresponds to $V_\delta (2)$, and the arrows label $V_\delta $, $\hat V_{1}$ and $\hat V_{2}$ at $t=22.5$. The second order approximation $\hat V_{2}$ holds on longer, but eventually succumbs. By then, $V_\delta $ is far from the constant Lie quadratic  $V_0$. $\square$
\end{example}
Theorems \ref{thm2}, \ref{thm3} give explicit formulae in terms of elementary functions for  approximations $\hat V_n$ to nearly constant Lie quadratics $V_\delta$ in $E^3$. Order $n$ approximations $x$ to Riemannian cubics $x_\delta$ in $SO(3)$ with Lie quadratic $V_\delta$ can be found by solving the linear system of ODEs
\begin{equation}\label{xquad}x^{(1)}(t)~=~x(t)\hat V_n(t)\end{equation}
where, on the right hand side, $\hat V_n(t)$ is identified with an element of $so(3)$ and, considered as a $3\times 3$ matrix, is premultiplied  by the unknown matrix $x(t)$. Because  
(\ref{xquad}) has nonconstant coefficients, $x$ cannot be written down directly and would be found by a numerical integrator such as NDSolve.  Considering that Riemannian cubics in $SO(3)$ are solutions of the $36$ dimensional nonlinear system (\ref{el}), the $9$ dimensional linear system (\ref{xquad}) appears comparatively benign. 
But solving (\ref{xquad}) is nontrivial, and this step can be avoided: an explicit first order approximation $\hat x$ to $x_\delta$ is given as follows. 
\section{Approximating Nearly Geodesic Cubics in $SO(3)$}\label{approxsecv}
Given a generic Lie quadratic in a semisimple Lie algebra, there is an integrability algorithm  \cite{lyleqjm}  that takes a single quadrature, and gives an explicit formula for an associated Riemannian cubic. The algorithm takes a simple form \cite{lyledual} for Riemannian cubics in $SO(3)$. So we might try to approximate a nearly geodesic 
Riemannian cubic $x_\delta :[t_0,t_1]\rightarrow SO(3)$ to first order by substituting $\hat V_1$ for $V_\delta $ in the integrability algorithm. This has the unexpected benefit of giving a first order approximation for the quadrature in terms of elementary functions. On the other hand, $\hat V_2$ is needed to approximate the other terms in the integrability algorithm to first order.  This gives, by somewhat indirect means, explicit first order approximations $\hat x$ to nearly geodesic Riemannian cubics $x$ in $SO(3)$, without even the need for a single quadrature. 
To review the integrability algorithm for Lie quadratics we make two definitions.  
\begin{definition} For $r\in \R$, let $R_r\in SO(3)$ be clockwise rotation by $r$ in $E^2$, namely  
$$R_r~:=~\left[ \begin{array}{ccc}
~\cos r&\sin r&0\\
-\sin r&\cos r&0\\
0&0&1
\end{array}\right] .$$
$\square$
\end{definition}
\begin{definition} For linearly independent $X_1,X_2\in E^3$, set 
$$S(X_1,X_2)~:=~\left[ \begin{array}{ccc}
\frac{\Vert X_1\Vert X_2-\langle X_1,X_2\rangle X_1 }{\sqrt{\Vert X_1\Vert ^2\Vert X_2\Vert ^2-\langle X_1,X_2\rangle ^2}}&\frac{X_1\times X_2}{\sqrt{\Vert X_1\Vert ^2\Vert X_2\Vert ^2-\langle X_1,X_2\rangle ^2}}&\frac{X_1}{\Vert X_1\Vert }
\end{array}\right] ^{\bf T}~\in ~SO(3)$$
where $~^{\bf T}~$ means transpose. Then, for any curve $W:[t_0,t_1]\rightarrow E^3$ with 
$W(t),W^{(1)}(t)$ everywhere linearly independent, define $T(W):[t_0,t_1]\rightarrow SO(3)$ by 
$$T(W)(t)~:=~S(W(t),W^{(1)}(t)).$$
$\square$
\end{definition}
Now we return to Lie quadratics\footnote{Here the symbol $V$ is used to denote a Lie quadratic, rather than a variation of Lie quadratics. } $V$. As is easily proved \cite{jmp},  they have constant acceleration: $c:=\langle V^{(2)}(t),V^{(2)}(t)\rangle $ is constant as $t$ varies. When $V^{(3)}$ vanishes identically in a nonempty open interval, the corresponding Riemannian cubics are cubically reparameterised geodesics. Assume $V^{(3)}(t)\not= {\bf 0}$ for any $t\in [t_0,t_1]$, where $V: [t_0,t_1]\rightarrow so(3)\cong E^3$, and $E^3$ is identified with $so(3)$ in the standard way, by the adjoint reprsentation. 
%
%
\begin{lemma}\label{lem1} Let $x:[t_0,t_1]\rightarrow SO(3)$ be a Riemannian cubic whose Lie quadratic is $V$. For $t\in [t_0,t]$, set   
$$y(t)~:=~R_{\phi (t)}\circ T(V^{(2)})(t)\quad \hbox{where}\quad 
\phi (t)~:=~c^{1/2}\int _{t_0}^t \frac{c-\langle C,V^{(2)}(s)\rangle }{\Vert V^{(3)}(s)\Vert ^2}~ds$$
and $\circ $ stands for matrix multiplication. Then $x(t)=x(t_0)y(t_0)^{\bf T}y(t)$. 
\end{lemma}

\spp
{\bf Proof:} Since left Lie reductions are invariant with respect to left multiplication, we can suppose without loss of generality that $x(t_0)$ is the identity. Setting $W_3(t):=c^{-1/2}V^{(2)}(t)$, 
$W_1(t):=V^{(3)}(t)/\Vert V^{(3)}(t)\Vert$, and $W_2(t):=W_3(t)\times W_1(t)$, we find that  $W_1,W_2,W_3$ meet the requirements of Theorem 5 of \cite{lyledual} (in \cite{lyledual} a precise choice of $W_1$ is not made, allowing other curves in $S^2$ orthogonal to $W_3$). So there is a Riemannian cubic $y$ of the form above, with Lie quadratic $V$, provided    
$$\phi ^{(1)}(t)~=~\langle W_1(t),W_2^{(1)}(t)+{\rm ad}^{-1}(V(t))\times W_2(t)\rangle ~ =~-\langle W_1^{(1)}(t)+{\rm ad}^{-1}(V(t))\times W_1(t),W_2(t)\rangle .$$   

\spp
With the present choice of $W_1$,  
~$\displaystyle{\phi ^{(1)}=\frac{1}{c^{1/2}\Vert V^{(3)}\Vert ^2}\langle V^{(3)},[V^{(2)}, V^{(4)}]+[V, [V^{(2)}, V^{(3)}]]\rangle =~}$ 
$$\frac{\langle V^{(3)},[V^{(2)}, [V^{(3)}, V]+[V^{(2)}, V^{(1)}]]+[V, [V^{(2)}, V^{(3)}]]\rangle}{c^{1/2}\Vert V^{(3)}\Vert ^2} ~=~\frac{\langle V^{(3)},[V^{(2)}, [V^{(2)},V^{(1)}]]\rangle}{c^{1/2}\Vert V^{(3)}\Vert ^2} ~=~$$
$$-c^{1/2}\frac{\langle V^{(3)},V^{(1)}\rangle }{\Vert V^{(3)}\Vert ^2}~=~c^{1/2}\frac{\langle V^{(2)},[V^{(1)}, V]\rangle }{\Vert V^{(3)}\Vert ^2}~=~c^{1/2}\frac{(c-\langle C,V^{(2)}\rangle )}{\Vert V^{(3)}\Vert ^2}.$$
$\square$

\spp
So Riemannian cubics in $SO(3)$ can be found from Lie quadratics by a single quadrature. First order approximations to nearly geodesic Riemannian cubics can be written explicitly in terms of elementary functions, without the need for quadrature. More precisely,  let $x_\delta :[t_0,t_1]\rightarrow SO(3)$ be a Riemannian cubic whose Lie quadratic $V_\delta :[t_0,t_1]\rightarrow so(3)$ is nearly constant. For $n=1,2$ let $\hat V_n:[t_0,t_1] \rightarrow so(3)$ be the order $n$ approximations to $V_\delta$ given in \S \ref{nearconstv}. Assuming $B\not= {\bf 0}$, write $B=\beta (\cos \gamma ,\sin \gamma )$ where $\gamma \in [0,2\pi )$.  Set $\rho :=-2c_2/(d^2\beta )$.
\begin{theorem}\label{thm4} For $t\in t_0,t_1]$ define  $\hat x(t):=x_\delta (t_0)\hat y(t_0)^{\bf T}\hat y(t)$ where $\hat y(t):=R_{\hat \phi (t)}\circ T(\hat V_2^{(2)})(t)$ and 
$$\hat \phi (t)~:=~ \delta (\rho ^2+1)^{1/2}((t-t_0)\beta +\frac{a_{11}(\cos (\gamma -d(t-t_0))-\cos \gamma )+a_{12}(\sin(\gamma -d(t-t_0))-\sin \gamma ) }{d^2 }).$$
Then for $j\geq 0$, and all $t\in [t_0,t_1]$, 
$x_\delta ^{(j)}(t)=\hat x^{(j)}(t)+O(\delta ^2)$, 
where the asymptotic constants depend on $a_{11},a_{12},\beta ,c_2,d,t_0,t_1$ and $j$.  
\end{theorem}

\spp
{\bf Proof:} Let $\hat W_3(t)$ and $\hat W_1(t)$ be the unit vectors in the directions of $\hat V_2^{(2)}(t)$ and $\hat V_2^{(3)}(t)$ respectively. Taking $V=V_\delta$ in Lemma \ref{lem1}, it suffices to show that $\phi (t)=\hat \phi (t)+O(\delta ^2)$ and that $W_i(t)=\hat W_i(t)+O(\delta ^2)$ for $i=3,1$. 

\spp
By Theorem \ref{thm2}, $V_\delta ^{(2)}=\hat V_1^{(2)}+O(\delta ^2)=\hat V_1^{(2)}(1+O(\delta ))$. So  
$c=\hat c(1+O(\delta ))$ and $C=\hat C(1+O(\delta ))$ where $\hat c:=\Vert \hat V_1^{(2)}\Vert ^2=\delta ^2(4c_2^2+d^4\beta  ^2)$, $\hat C:=\delta ~{\rm ad}(2c_2,-da_{12},da_{11})=\hat V_1^{(2)}-[\hat V_1^{(1)},\hat V_1]+O(\delta ^2)$. 
Similarly $V_\delta ^{(3)}=\hat V_1^{(3)}(1+O(\delta ))\Longrightarrow \Vert V_\delta ^{(3)}\Vert ^2=d^6\beta ^2+O(\delta )$.  So 
$$\phi (t)~=~ c^{1/2}\int _{t_0}^t \frac{c-\langle C, V_\delta ^{(2)}(s)\rangle }{\Vert V_\delta ^{(3)}(s)\Vert ^2}~ds~=~\hat c^{1/2}\int _{t_0}^t \frac{\hat c-\langle \hat C, \hat V_1^{(2)}(s)\rangle }{\Vert \hat V_1^{(3)}(s)\Vert ^2}~ds +O(\delta ^2).$$
Because the denominator is constant, the integral on the right can be computed precisely, giving 
\begin{eqnarray*}\phi (t)&=&\delta (4c_2^2+d^4\beta  ^2)^{1/2}(\frac{t-t_0}{d^2}+\frac{a_{11}(\cos (\gamma -d(t-t_0))-\cos \gamma )+a_{12}(\sin(\gamma -d(t-t_0))-\sin \gamma ) }{d^4\beta })+O(\delta ^2)\\
 &=&\hat \phi (t)+O(\delta ^2).\end{eqnarray*}
 
 \spp
 By Theorem \ref{thm3}, $V_\delta ^{(2)}=\hat V_2^{(2)}+O(\delta ^3)=\hat V_2^{(2)}(1+O(\delta ^2))$. So $W_3=\hat W_3+O(\delta ^2)$. Similarly $W_1=\hat W_1+O(\delta ^2)$. $\square$

 \spp
Whereas in \cite{lylenonnull} the analysis of Riemannian cubics in $SO(3)$ is 
complicated,  in the present paper $\hat x$ is algebraic in terms of polynomials and trigonometric functions, as we see by taking $n=2$ in equations (\ref{forfn}), (\ref{forvn}):
\begin{eqnarray*}
f_2^{(2)}&=&2d[A_0,\tilde i(\tilde e-{\bf 1})B]+[A_1,2(\tilde e-{\bf 1}+d(t-t_0)\tilde i \tilde e)B] \\
v_2^{(2)}&=&-\frac{4c_2}{d}(\tilde e-{\bf 1})A_0+\frac{4c_2}{d^2}(d(t-t_0){\bf 1}-\tilde i(\tilde e-{\bf 1}))A_1+2(2c_2(t-t_0)+d^2I(q))\tilde i\tilde eB 
%
\end{eqnarray*}
where everything is constant, except  $t$,
$$\tilde e(t)~=~ \left[ \begin{array}{ccc}
1&0&0\\
0&~\cos (d(t-t_0))&\sin (d(t-t_0))\\
0&-\sin (d(t-t_0))&\cos (d(t-t_0))
\end{array}\right] \quad \hbox{and}\quad I(q)(t)~=~\frac{c_2(t-t_0)^3}{3}+\frac{c_1(t-t_0)^2}{2}+c_0(t-t_0).$$
Then   
$\hat V_2^{(2)}=\delta (2c_2F_0-d^2\tilde eB)+\frac{\delta ^2}{2}(f_2^{(2)}F_0
+v_2^{(2)})
$ is affine in $\tilde e(t)B$, $f_2^{(2)}(t)$ and $v_2^{(2)}(t)$. So in Theorem \ref{thm4}, $T(\hat V_2^{(2)})$ is $S(X_1,X_2)$ where 
\begin{eqnarray*}
X_1(t)&:=&2c_2F_0-d^2\tilde e(t)B+\frac{\delta}{2}(f_2^{(2)}(t)F_0+v_2^{(2)}(t))\quad \hbox{and}\\
X_2(t)&:=&d^3\tilde i\tilde e(t)B+\frac{\delta }{2}(f_2^{(3)}(t)+v_2^{(3)}(t)).
\end{eqnarray*}
\begin{example}\label{ex2} For some small $\delta \in \R$, let $V_\delta :[0,10]\rightarrow so(3)$ be the Lie quadratic satisfying 
\begin{equation}\label{icsv}{\rm ad}(V_{\delta}(0))~=~(1,0,0)+\delta (0,1,0) ,\quad {\rm ad}(V_\delta ^{(1)}(0))~=~\delta (0,0,1)/2 ,\quad {\rm ad}(V_\delta {(2)}(0))~=~\delta (1,1,1)/4 \end{equation}
and let $x_\delta :[0,8]\rightarrow SO(3)$ be the corresponding Riemannian cubic for which $x_\delta (0)={\bf 1}$. 

\spp
Given any particular value of $\delta$, say $\delta =0.05$, Mathematica's NDSolve can be used to numerically solve 
the quadratic differential equation (\ref{el}) for $V_\delta :[0,8]\rightarrow so(3)$. Then $x_\delta $ is found by numerically solving the linear differential equation 
$$x_\delta ^{(1)}(t)~=~x_\delta (t)V_\delta (t)$$
where the coefficients on the right are entries of the matrix $V_\delta (t)$.  The second rows of $x_\delta (t)$ are shown as  
%
%
 the blue curve  in 
Figure \ref{fig:approx_cubic}, with blue labelled points corresponding to $t=0,1,2,\ldots 10$.  A geodesic represented in this way would appear as a circular arc. 

\spp
To compute the approximate solution $\hat x:[0,10]\rightarrow SO(3)$ from Theorem \ref{thm4}, take $t_0=0,t_1=10$ and $D=(1,0,0)$. For $\hat x$ with $\hat x(0)={\bf 1}$ to satisfy the other  initial conditions (\ref{icsv}) with $O(\delta ^2)$ errors, it suffices to have   
$$\hat V_1(0))~=~(1,\delta ,0) ,\quad \hat V_1^{(1)}(0)~=~(0,0,\delta /2) ,\quad \hat V_1{(2)}(0)~=~\delta (1,1,1)/4 $$
where $\hat V_1$ is considered as a curve in $E^3$.  For this, choose parameters $c_0,c_1,c_2,a_{01},a_{02},a_{11},a_{12},\beta ,\gamma $ so that 
\begin{eqnarray*}
(c_0,a_{01}+\beta \cos \gamma ,a_{02}+\beta \sin \gamma )&=&(0,1 ,0)\\
(c_1,a_{11}+\beta \sin \gamma ,a_{12}-\beta \cos \gamma )&=&
(0,0,1)/2\\
(c_2,-\beta \cos \gamma ,-\beta \sin \gamma )&=&(1,1,1)/4 
\end{eqnarray*}
Taking  $\beta =\sqrt{2}/4$ we have $c_0=c_1=0$, $c_2=1/8$, $\gamma =5\pi /4$, $a_{01}=5/4$, $a_{02}=1/4$, 
$a_{11}=a_{12}=1/4$. 
\begin{figure}[p] 
   \centering
   \includegraphics[width=3.5in]{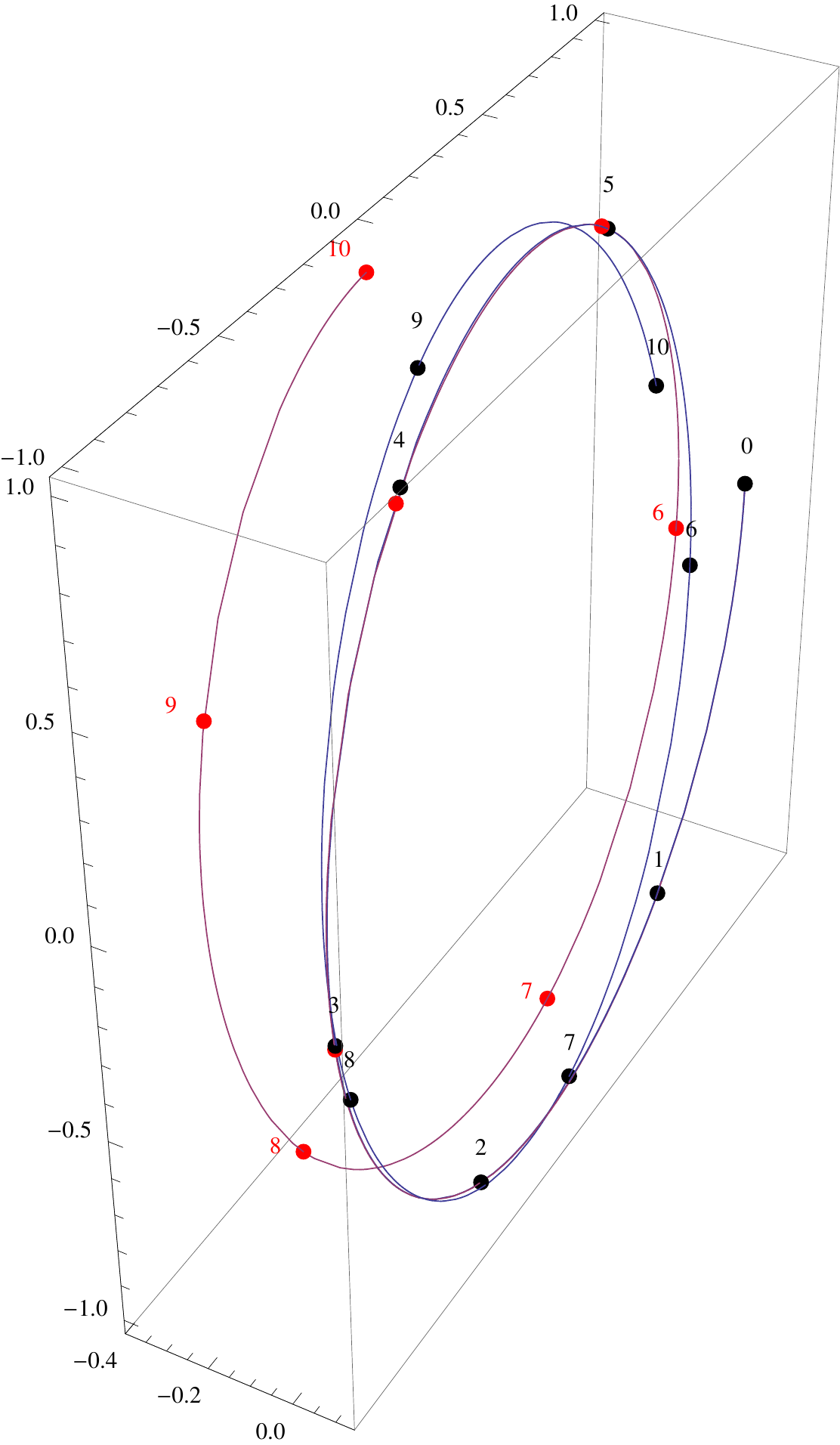} 
   \caption{{\rm Approximation (second rows)  to a 
   Nearly Geodesic Riemannian Cubic in} $SO(3)$}
   \label{fig:approx_cubic}
\end{figure}
The second rows of $\hat x(t)$ are shown as the red curve in Figure \ref{fig:approx_cubic}, with red points corresponding to $t=3,4,5,\ldots 10$. For $t=0,1,2$, $\hat x(t)$ and $x_\delta (t)$ are hard to distinguish. For $t=3,4,5$ the approximation is good enough for the blue numbers to label both.

\spp
It is hard to see much difference between $\hat x(t)$ and $x_\delta (t)$ until around $t=4$, and long before then the Riemannian cubic $x_\delta$ is obviously nongeodesic. After that, differences between $\hat x(t)$ and $x_\delta (t)$ become noticeable ($t=5$), and then large ($t\geq 6$) where  $x_\delta$ is very far from geodesic. 
$\square$
\end{example}  
\section{Conclusion}
We derive first order approximations $\hat x$ to nearly geodesic Riemannian cubics $x_\delta $ in $SO(3)$ defined over a given interval $[t_0,t_1]$. The approximations capture a good deal of the geometry of $x_\delta$, including behaviour of the associated derivative, corresponding to the body angular velocity for a rigid body. Equivalently, the approximations apply to restrictions of arbitrary Lie quadratics to  sufficiently small intervals. This fills a gap in the literature, by complementing studies on long term asymptotics of Riemannian cubics. 

\spp


\begin{thebibliography}{9}
\bibitem{oscar}
J. Arroyo, O.J. Garay and J. Mencia, 
``Elastic circles in $2$-spheres," 
{\em J. Phys. A: Math. Gen.} {\bf 39} (2006) 2307--2324.

\bibitem{bloch}
A.M. Bloch with J. Baillieul, P. Crouch and J. Marsden, 
{\em Nonholonomic Mechanics and Control}, 
Interdisciplinary Applied Mathmatics, Springer (2003).

\bibitem{crouch1} 
M. Camarinha, F. Silva Leite, P. Crouch,
``On the geometry of Riemannian cubic polynomials,'' 
{\em Differential Geom. Appl.} {\bf 15} (2) (2001) 107--135.

\bibitem{crouch4}
P. Crouch and F. Silva Leite,  
``The dynamic interpolation problem: on Riemannian manifolds, Lie groups, 
and symmetric spaces,''
{\em J. Dynam. Control Systems} {\bf 1} (2) (1995) 177--202. 



\bibitem{giambo}
R. Giambo, F. Giannoni and P. Piccione, 
``An analytical theory for Riemannian cubic polynomials," 
{\em IMA J. Math. Control \& Information} {\bf 19} (2002) 445--460.

\bibitem{barr}
R. Ramamoorthi and A. Barr, 
``Fast construction of accurate quaternion splines," 
{\em Proc. of SIGGRAPH} '97, Los Angeles, August 3Ð8, (1997) 287--292. 

\bibitem{vega}
S. Guti\'errez, J. Rivas and L. Vega,
``Formation of singularities and self-similar vortex motion under the localized induction approximation,"
{\em Comm. in Partial Differential Equations} {\bf 28} (2003) 927--968. 


\bibitem{jurd}
V. Jurdjevic,
``Non-Euclidean elastica," 
{\em Amer. J. of Math.} {\bf 117} (1995) 93--124.

\bibitem{lylegreg}
L. Noakes, G. Heinzinger and B. Paden,
``Cubic splines on curved spaces,''
{\em IMA J. Math. Control \& Information} {\bf 6} (1989) 465--473.

\bibitem{jmp}
L. Noakes,
``Null cubics and Lie quadratics,''
{\em J. Math. Physics} {\bf 44} (3) (2003) 1436--1448. 

\bibitem{lylenonnull}
L. Noakes,
``Non-null Lie quadratics in $E^3$,'' 
{\em J. Math Physics}, {\bf 45} (11) (2004) 4334--4351.

\bibitem{lyledual}
L. Noakes,
``Duality and Riemannian Cubics,'' 
{\em Adv. in Computational Math.} {\bf 25} (2006) 195--209.

\bibitem{lyleqjm}
L. Noakes,
``Lax constraints in semisimple Lie groups,''
{\em Quart. J. Math.}, {\bf 57} (2006) 527--538. 

\bibitem{lyleSIAM}
L. Noakes, 
``Asymptotics of null Lie quadratics in $E^3$,'' 
{\em SIAM J. on Applied Dynamical Systems}, {\bf 7} (2) (2008) 437--460. 

\bibitem{lyletomaszrev}
L. Noakes and T. Popiel,
``Geometry for robot path planning," 
{\em Robotica} {\bf 25} (2007) 691--701.

\bibitem{pauley}
M. Pauley,   
``Null Lie quadratics in $so(3)$ and $sl(2)$: applications and quadratures," 
{\em in preparation}, (February 2009). 

\bibitem{elas}
T. Popiel and L. Noakes, 
``Elastica in $SO(3)$," 
{\em J. Australian Math. Soc.} {\bf 83} (2007) 105--125.

\bibitem{singer}  
D.A. Singer,
``Lectures on elastic curves and rods,'' 
{\em Curvature and Variational Modeling in Physics and Biophysics,  AIP Conf. Proc.} {\bf 1002} (2008) 3--32. 
%

\bibitem{vara}
V.S. Varadarajan,
{\em Lie Groups, Lie Algebras and Their Representations},
Prentice-Hall 1974. 
%
\end{thebibliography}
\end{document}